%
%
\documentstyle{amsppt}
\topmatter
\parskip=0.5\baselineskip 

\loadmsbm
\UseAMSsymbols
\hoffset=0.75truein
\voffset=0.5truein

\def\=#1{\accent"16 #1}
\def\x{\relax}
\def\forget#1{}

\def\deltaeff{\operatorname{\delta_{\operatorname{\scriptscriptstyle eff}}}}
\def\doct{\delta_{oct}}
\def\pt{\hbox{\it pt}}
\def\Vol{\hbox{vol}}

\def\sqr{\sqrt}
\def\sol{\operatorname{sol}}
\def\vor{\operatorname{vor}}
\def\dih{\operatorname{dih}}

\def\cro{\operatorname{crown}}
\def\anc{\operatorname{anc}}
\def\quo{\operatorname{quo}}

\def\rad{\operatorname{rad}}

\def\diag|#1|#2|{\vbox to #1in {\vskip.3in\centerline{\tt Diagram #2}\vss} }
\def\v{\hskip -3.5pt }
\def\gram|#1|#2|#3|{
        {
        \smallskip
        \hbox to \hsize
        {\hfill
        \vrule \vbox{ \hrule \vskip 6pt \centerline{\it Diagram #2}
         \vskip #1in %
             \includegraphics{#3}\hrule }
        \v\vrule\hfill
        }
\smallskip}}

\title A Formulation of the Kepler Conjecture\endtitle

\author Samuel P. Ferguson, Thomas C. Hales\endauthor
\endtopmatter
\document

The Kepler conjecture asserts that the density of a packing
of equal spheres cannot be greater than that of the face-centered
cubic packing.   There are various optimization problems in
finitely many variables that imply the Kepler conjecture.  The
first was introduced by L. Fejes T\'oth \cite{FT}.  His formulation
was based on the Voronoi decomposition of space.  In \cite{H1}, a dual
formulation based on the Delaunay decomposition was proposed.
Later, in \cite{I}, the two strategies were combined to
partition space by a combination of the Voronoi and Delaunay decompositions.
\footnote""{\hfill\it version - 7/18/98, corrected -12/5/01}
\footnote""{\noindent Research supported in part by the NSF}

The Voronoi and Delaunay decompositions can be mixed
in infinitely many ways.  This gives a large set
of optimization problems in finitely many variables
that imply the Kepler conjecture.  Each of these optimization
problems is a formulation of the Kepler conjecture.

The selection of a good formulation of the Kepler conjecture
is a central issue in the resolution of the conjecture.
Our experience suggests that
whenever the technical difficulties become too great it is
generally better to rework the formulation of the problem
rather than to confront the technical difficulties directly.
It is the infinite dimensionality of the problem that
gives the flexibility to skirt these technical problems.

In \cite{I}, five steps were suggested, which collectively imply
the Kepler conjecture.  The formulation suggested in \cite{I}
was sufficient to see the first two steps to completion.
To see the third and fifth steps to completion, it has been
necessary to make some adjustments.
This paper makes those changes.  The proofs of the third and
fifth steps, found in \cite{III} and \cite{V}, rely essentially
on the constructions and results of this paper.  We have been
careful to modify the constructions in a way that does not
affect the proofs of the first two steps of the program.  Lemma
3.13, Proposition 3.14, Conjecture 3.15, and Theorem 3.16 are the results
needed to bring the results of \cite{I} and \cite{II} into
harmony with this paper.  Proposition 4.1 is used in \cite{III}
and \cite{V}.  Proposition 4.7 brings significant simplifications
to the calculations of \cite{III}.

Our formulation has departed more than ever from the original
formulation on the space of Delaunay stars (defined in \cite{H1}).
The Delaunay decomposition plays a smaller role here than in
any of our previous papers, although many of the concepts it
inspired remain (such as the compression of a simplex,
quasi-regular tetrahedra, and quarters).
To reflect this change in formulation, we now call the
stars {\it decomposition stars}.

\head 1. Geometric Decomposition\endhead

Fix a packing of spheres of radius 1.
The centers of
the spheres are called  {\it vertices}.  Fix the constant
$2.51$.  It is used throughout this paper and all of our
related papers on the subject.  A {\it quasi-regular\/}
tetrahedron is the simplex formed by four vertices,
 each at most $2.51$ from
the others.
A {\it quarter\/} is defined as a
simplex whose edge lengths $y_1,\ldots,y_6$ can be ordered to
satisfy $2.51\le y_1\le \sqr8$, $2\le y_i\le 2.51$, $i=2,\ldots,5$.
We call the longest
edge of a quarter its {\it diagonal\/}.
When the quarter has a distinguished vertex, the
quarter is {\it upright\/} if the distinguished vertex is an endpoint
of the diagonal, and {\it flat\/} otherwise.

If four quarters fit together along a common diagonal, forming a figure
with six vertices, the
resulting figure is called an {\it octahedron}.
The
octahedron may have more than one diagonal of length at most $\sqr8$,
so its decomposition into four quarters need not be unique.
(This definition of octahedron differs from the one given in \cite{I},
but we will never return to the earlier definition.)

Our simplices are generally assumed to come with a distinguished
vertex, fixed  at the origin.
(The origin will always be at a vertex of the packing.)
We
number the edges of each simplex $1,\ldots,6$, so that
edges $1$, $2$, and $3$ meet at the origin, and
the edges $i$ and $i+3$ are opposite, for $i=1,2,3$.
$S(y_1,y_2,\ldots,y_6)$ denotes a simplex whose edges have
lengths $y_i$, indexed in this way.
We refer to the endpoints away from the origin
of the first, second, and third
edges as the first, second, and third
vertices.

We say that two manifolds with boundary {\it overlap\/} if their
interiors intersect. We define the projection of a set $X$ to be
the radial projection of $X\setminus{0}$ to the unit sphere
centered at the origin. We say they {\it cross\/} if their
projections to the unit sphere overlap.  We label the edges of a
simplex $S(y_1,\ldots,y_6)$ as in \cite{I}. In general, let
$\dih(S)$ be the dihedral angle of a simplex $S$ along its first
edge. When we write a simplex in terms of its vertices
$(w_1,\ldots,w_4)$, then $(w_1,w_2)$ is understood to be the first
edge. We define a function $\Cal
E(S(y_1,\ldots,y_6),y_1',y_2',y_3')$, by taking two simplices
$S=S(y_1,\ldots,y_6)$ and $S'=S(y_1',y_2',y_3',y_4,y_5,y_6)$,
 and moving $S'$ until the simplices do not overlap,
and the face formed by the fourth, fifth, and sixth edges of
$S$ and $S'$ coincide.
$\Cal E$ is defined only if $S$ and $S'$ exist, and then it
is defined as the distance between
the origin and the vertex $v'$ of $S'$ opposite the common face
(Diagram 1.1).

\smallskip
\gram|1.5|1.1|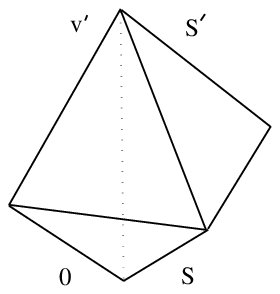|
\smallskip

If intervals containing $y_1,\ldots,y_6,y_1',y_2',y_3'$ are given,
 lower
bounds on $\Cal E$ over that domain are generally easy to obtain.
For example,
if the segment from the vertex $v'$ of $S'$ to the origin passes
through the face common to $S$ and $S'$, then
$\Cal E$ is
increasing in the variables $y_1,y_2,y_3,y_1',y_2',y_3'$ (at least until we
deform the simplices sufficiently that the segment no longer
passes through the common face).  A {\it pivot\/} is the circular motion
of a vertex at a fixed distance from two others  (see \cite{I}).
The {\it axis\/} of the pivot is the line through the two fixed vertices.
By using pivots,
we observe that $\Cal E$ is monotonic
decreasing in $y_4,y_5,y_6$.
For example, if we pivot the first vertex away from the third
around the axis through the second edge, $\Cal E$ is unaffected.
Because these lower bounds are generally so easily determined,
we will state them without proof.  We will state that these
bounds were obtained by {\it geometric considerations}, to
indicate that the bounds were obtained by the monotonicity arguments
of this paragraph.

\bigskip
\proclaim{Lemma 1.2}  No vertex of the packing is contained in the
interior of a quasi-regular tetrahedron or quarter.
\endproclaim

\demo{Proof} See I.3.5. \qed\enddemo

\proclaim{Corollary}  No vertex of the packing is contained in the
interior of an octahedron. \qed\endproclaim

\proclaim{Lemma 1.3}  An edge of length $2.51$ or less cannot pass
through a face whose edges have lengths $2.51$, $2.51$, and $\sqr8$ or
less.
\endproclaim

\demo{Proof}  The distance between each pair of vertices is at least 2.
Geometric considerations
show that the edge has length at least
$${\Cal E}(S(2,2,2,2.51,2.51,\sqr8),2,2,2) > 2.51.$$
\qed\enddemo

\proclaim{Lemma 1.4}  If the diagonal of a quarter passes through a face
of a quasi-regular tetrahedron, then each of the two endpoints of the
diagonal edge is at most 2.2 away from each of the vertices of the face
(see Diagram I.3.1).
\endproclaim

\demo{Proof}  Let the diagonal edge be $(w_1,w_2)$ and the vertices of
the face be $(v_1,v_2,v_3)$.  If $|v_i-w_j|>2.2$, for some $i$, $j$,
then geometric considerations give
$$|w_1-w_2|\ge{\Cal E}(S(2,2,2,2.51,2.51,2.51),2,2,2.2) > \sqr8.$$\qed
\enddemo

As in earlier papers,
$\eta(x,y,z)$ denotes the circumradius of a triangle with
edge-lengths $x$, $y$, and $z$.
Suppose that $S$ and $S'$ are adjacent quasi-regular tetrahedra
with a common face $F$.  As in Lemma 1.4, suppose that a diagonal
of a quarter runs between the opposite vertices of $S$ and $S'$
through the face $F$.
By the lemma, each of the six external faces of the
pair of quasi-regular tetrahedra
has circumradius at most $\eta(2.2,2.2,2.51)<\sqr2$.
A diagonal of a quarter  cannot pass through a face of this size
\cite{I.3.2}.
This pair of quasi-regular tetrahedra is the union of three quarters
joined along a common diagonal.   No other quarter overlaps these
quasi-regular tetrahedra.

If $(v_1,v_2)$ is an edge of length between 2.51 and $\sqr8$,
we say that $v$ $(\ne v_1,v_2)$ is an {\it anchor\/} of $(v_1,v_2)$ if
its distances to $v_1$ and $v_2$ are at most $2.51$.
The two vertices of a quarter that
are not on the diagonal are anchors of the diagonal,
and the diagonal may have other
anchors as well.

\proclaim{Lemma 1.5}  Suppose an edge $(w_1,w_2)$ of length at most $\sqr8$
passes through the face formed by a diagonal $(0,v_1)$ and one of its anchors.
Then $w_1$ and $w_2$
are also anchors of $(0,v_1)$.
\endproclaim

\demo{Proof}  ${\Cal E}(S(2,2,2,\sqr8,2.51,2.51),2,2,2.51) > \sqr8$.
\qed\enddemo

The {\it height\/} of a vertex is its distance
from the origin.
We say that a vertex is {\it enclosed\/} over a figure
if it lies in the interior of the
cone at the origin generated by the figure.

  If we draw a geodesic arc on the unit sphere with endpoints at the
projections of $v_1$ and $v_2$
for every pair of vertices $v_1$, $v_2$ such that
$|v_1|, |v_2|, |v_1-v_2|\le 2.51$, we obtain a planar map that breaks
the unit sphere into regions called {\it standard regions}.
(The arcs do not meet except at endpoints \cite{I.3.10}.)

By a pair of adjacent quarters, we mean two quarters
sharing a face along the diagonal.
The common vertex that
does not lie on the diagonal is called the {\it base point\/} of
the pair of adjacent quarters.  The other four vertices
are called the {\it corners\/} of the configuration.

\proclaim{Lemma 1.6}
Suppose that there exist four vertices $v_1,\ldots,v_4$ of height
at most 2.51 (that is, $|v_i|\le 2.51$) forming a skew quadrilateral.
Suppose that the diagonals $(v_1,v_3)$ and $(v_2,v_4)$
have lengths between $2.51$ and $\sqr8$.
Suppose the diagonals $(v_1,v_3)$ and $(v_2,v_4)$
cross.
Then the four vertices are the corners of a pair of adjacent
quarters with base point at the origin.
\endproclaim

\demo{Proof}
 Set $d_1=|v_1-v_3|$ and $d_2 = |v_2-v_4|$.
By hypothesis, $d_1$ and $d_2$ are at most $\sqr8$.
 If $|v_1-v_2|>2.51$,
geometric considerations give the contradiction
$$\max(d_1,d_2)\ge \Cal E(S(2.51,2,2,2.51,\sqr8,2.51),2,2,2) > \sqr8
    \ge\max(d_1,d_2).$$
Thus, $(0,v_1,v_2)$ determines a bounding arc of
standard region,  as do $(0,v_2,v_3)$, $(0,v_3,v_4)$, and $(0,v_4,v_1)$
by symmetry. \qed\enddemo

\proclaim{Lemma 1.7}  If, in the context of Lemma 1.6, there is a
vertex $w$ of height at most $\sqr8$ enclosed over the
pair of adjacent quarters,
then $(0,v_1,\ldots,v_4,w)$ is an octahedron.
\endproclaim

\demo{Proof}  If the enclosed $w$ lies over
say $(0,v_1,v_2,v_3)$,
then $|w-v_1|$, $|w-v_3|\le 2.51$ (Lemma 1.5), where $(v_1,v_3)$ is
a diagonal.  Similarly, the distance from $w$ to the other two
corners is at most $2.51$.\qed\enddemo

\bigskip

We will select a nonoverlapping collection
of quarters and quasi-regular tetrahedra,  called
a {\it $Q$-system\/} (for quarters and quasi-regulars).
For each octahedron, we fix a diagonal
of length at most $\sqr8$ and
place the four quarters along that diagonal in the $Q$-system, but not
the overlapping quarters situated
along other diagonals of the octahedron.
Place all quasi-regular tetrahedra in the $Q$-system. This, of course,
prevents us from placing
any quarters that overlap these tetrahedra into the
$Q$-system (as in Lemma 1.4).

Fix the origin at the base point of a pair of adjacent quarters.
We investigate the local geometry when another quarter overlaps
one of them.  This happens, for example, if both diagonals between
opposite corners of the pair of quarters have lengths
at most $\sqr8$.  We will see that a conflict like this between the
diagonals between corners is the
only way a pair of adjacent quarters can overlap another quarter.
We call these conflicting diagonals.
Label the four corners of the pair of quarters
$v_1$, $v_2$, $v_3$, $v_4$,
with $(v_1,v_3)$ the common diagonal.
We say that an edge is {\it short\/}
if its length is at most $2.51$.

\noindent{\bf Case 1.}
There is an enclosed vertex $w$, say over
$(v_1,v_2,v_3)$, where $(0,w)$ is a diagonal of a quarter.
Lemma 1.5 implies that $v_1$ and $v_3$ are anchors of $(0,w)$.
The only other possible anchors of $(0,w)$ are $v_2$ or $v_4$,
for otherwise a short edge passes through a face formed by $(0,w)$ and
one of its anchors.
If both $v_2$ and $v_4$ are anchors, then we have an
octahedron.  Otherwise, $(0,w)$ has at most $3$ anchors: $v_1$, $v_3$,
and either $v_2$ or $v_4$.  In fact, it must have exactly three
anchors, for otherwise there is no quarter along the edge $(0,w)$.
So there are exactly two quarters along the edge $(0,w)$.
We place the quarters along the diagonal $(v_1,v_3)$
in the $Q$-system.
The other two quarters, along the diagonal $(0,w)$, are not placed
in the $Q$-system.  They form a pair of adjacent quarters
(with base point $v_4$ or $v_2$) that has conflicting diagonals,
$(0,w)$ and $(v_1,v_3)$,
of length at most $\sqr8$.

\noindent
{\bf Case 2.}  $(v_2,v_4)$ is a diagonal of
length at most $\sqr8$ (conflicting with $(v_1,v_3)$).  By
symmetry, we may assume that $(v_2,v_4)$ passes through the face $(0,v_1,v_3)$.
Assume (for a contradiction) that both diagonals have an
anchor other than the corners $v_i$.
Let the anchor of $(v_2,v_4)$ be denoted
$v_{24}$ and that of $(v_1,v_3)$ be $v_{13}$.
Assume the figure is not an octahedron, so that
$v_{13}\ne v_{24}$.
By Lemma 1.3, it is impossible to draw the edges $(v_1,v_{13})$
and $(v_{13},v_3)$ between $v_1$ and $v_3$.  In fact, if the
edges pass outside the quadrilateral $(0,v_2,v_{24},v_4)$, one of the
short edges $(0,v_2)$, $(v_2,v_{24})$, $(v_{24},v_4)$, or $(v_4,0)$
violates the lemma applied to the face $(v_1,v_3,v_{13})$.  If they
pass inside the quadrilateral, one of the edges $(v_1,v_{13})$, $(v_{13},v_3)$
violates the lemma applied to the face $(0,v_{2},v_4)$ or $(v_{24},v_2,v_4)$.
 We conclude that at most one of the two diagonals
has additional anchors.

If neither of the two diagonals has more than three anchors, we
have nothing more than two overlapping pairs of adjacent quarters along
conflicting diagonals.  Place the
two quarters along the lower edge $(v_2,v_4)$ into
the $Q$-system.
If there is a diagonal with more than three anchors,
place the
quarters along the diagonal with more than three anchors in the
$Q$-system.
  In both possibilities of case 2, the two quarters
left out of the $Q$-system correspond to a conflicting diagonal.

By the following lemma, Cases 1, 2, and the octahedron
are the only possibilities for a pair of adjacent quarters.

\proclaim{Lemma 1.8}
Let $v_1$ and $v_2$ be anchors of $(0,w)$ with $2.51\le |w|\le \sqr8$.
If an edge $(v_3,v_4)$
passes through both faces, $(0,w,v_1)$ and
$(0,w,v_2)$, then $|v_3-v_4|>\sqr8$.
\endproclaim

\demo{Proof}  Suppose the figure exists with $|v_3-v_4|\le\sqr8$.
Label vertices so $v_3$ lies on the same side of the
figure as $v_1$.
Contract $(v_3,v_4)$ by moving $v_3$ and $v_4$ until
    $(v_i,u)$ has length $2$,
for $u=0,w,v_{i-2}$, and $i=3,4$.
Pivot $w$ away from $v_3$ and $v_4$ around the axis $(v_1,v_2)$ until
    $|w|=\sqr8$.
Contract $(v_3,v_4)$ again. By stretching $(v_1,v_2)$,
we obtain a square of edge two and vertices
$(0,v_3,w,v_4)$.
Short calculations based on I.8.3.1 and its partial derivatives give
$$\dih(S(\sqr8,2,y_3,2,y_5,2)) > 1.075,\quad y_3,y_5\in[2,2.51],\tag1.7.1$$
$$\dih(S(\sqr8,y_2,y_3,2,y_5,y_6)) >1,\quad
        y_2,y_3,y_5,y_6\in[2,2.51].\tag1.7.2$$
Then
$$\pi\ge \dih(0,w,v_3,v_1) + \dih(0,w,v_1,v_2) + \dih(0,w,v_2,v_4)
    > 1.075 + 1 + 1.075 > \pi.$$
Therefore, the figure does not exist.\qed\enddemo

\proclaim{Lemma 1.9}  Let $v_1,v_2,v_3$ be anchors of $(0,w)$, where
$2.51\le|w|\le\sqr8$, $|v_1-v_3|\le\sqr8$, and the edge $(v_1,v_3)$
passes through the face
$(0,w,v_2)$.  Then $\min(|v_1-v_2|,|v_2-v_3|)\le2.51$.
Furthermore, if the minimum is $2.51$, then $|v_1-v_2|=|v_2-v_3|=2.51$.
\endproclaim

\demo{Proof}
Assume $\min\ge2.51$.
As in the proof of Lemma 1.8, we may assume that $(0,v_1,w,v_3)$ is
a square.  We may also
assume, without loss of generality, that $|w-v_2|=|v_2|=2.51$.
This forces $|v_2-v_i|=2.51$, for $i=1,3$.
This is rigid,  and is the unique figure that satisfies the constraints.
The lemma follows.\qed\enddemo

Assume that there are two quarters $Q_1$ and $Q_2$ that overlap.
Assume that neither is adjacent to another quarter.
Let $(0,u)$ and
$(v_1,v_2)$ be the diagonals of $Q_1$ and $Q_2$.
Suppose the diagonal $(v_1,v_2)$ passes through a face $(0,u,w)$ of $Q_1$.
By Lemma 1.5, $v_1$ and $v_2$ are anchors
of $(0,u)$.  Again, either the length of $(v_1,w)$ is at most
2.51 or the length of $(v_2,w)$ is at most 2.51,
say $(w,v_2)$.
It follows that $Q_1=(0,u,w,v_2)$ and $|v_1-w|\ge2.51$.
($Q_1$ is not adjacent
to another quarter.)  So $w$ is not an anchor of $(v_1,v_2)$.

Let $(v_1,v_2,w')$ be a face of $Q_2$ with $w'\ne 0,u$.
If $(v_1,w',v_2)$ does not link $(0,u,w)$, then $(v_1,w')$ or
$(v_2,w')$ passes through the face $(0,u,w)$, which is impossible by
Lemma 1.3.
So $(v_1,v_2,w')$  links $(0,u,w)$ and an edge of $(0,u,w)$ passes through
the face $(v_1,v_2,w')$.
It is not the edge $(u,w)$ or $(0,w)$, for they
are too short by Lemma 1.3.  So $(0,u)$ passes through $(w',v_1,v_2)$.
The only other anchors of $(v_1,v_2)$ are $u$ and $0$ (by Lemma 1.8).
Either $(u,w')$ or $(w',0)$ has length at most $2.51$ by Lemma 1.9, but
not both,
because this would create a quarter adjacent to $Q_2$.
By symmetry,
$Q_2=(v_1,v_2,w',0)$ and the length of $(u,w')$ is greater than 2.51.
By symmetry, $(0,u)$ has no other anchors either.
This determines the local geometry when there
are two quarters that intersect without belonging to a pair
of adjacent quarters
(see Diagram 1.10).

\smallskip
\gram|2|1.10|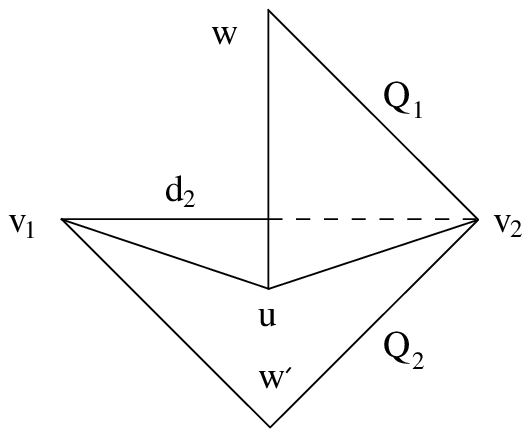|
\smallskip

When there are two isolated but overlapping quarters $Q_1$
and $Q_2$, then we place neither in the $Q$-system.  (We call
such a configuration an {\it isolated pair}.)
This completes the specification of the $Q$-system.  By construction,
if any quarter along a diagonal lies in the $Q$-system, then
all quarters along the diagonal lie in the $Q$-system.

\proclaim{Lemma 1.11}  Two vertices of height at most $\sqr8$ cannot be
enclosed over a flat quarter.\endproclaim

\demo{Proof} Assume the figure exists.
The diagonal $(v_1,v_2)$ of the quarter $(0,v_1,v_2,v_3)$ has
anchors $(0,v_3,w,w')$.  Lemma 1.8 gives $|w'|>\sqrt8$.
\qed\enddemo
%

\head 2. Voronoi Cells\endhead

In this section,  we
show that a mild modification of the Voronoi cells, called $V$-cells,
is compatible with the $Q$-system.

Recall from Section I.8.2, that the {\it orientation\/} of the face of a simplex
is said to be negative if the plane through that face separates the
circumcenter of the simplex from the vertex of the simplex that does
not lie on the face.  The orientation is positive if the circumcenter
and the vertex lie on the same side of the plane.

\proclaim{Lemma 2.1}   At most one face of a quarter $Q$ has
negative orientation.
\endproclaim

\demo{Proof}  The proof applies to any simplex with nonobtuse faces.
 Fix an edge and project $Q$ to a triangle in a
plane perpendicular to that edge.
The faces $F_1$ and $F_2$ of $Q$ along the edge
project to edges $e_1$ and $e_2$ of the triangular projection of $Q$.
The line equidistant from the three vertices of $F_i$
projects to a line perpendicular to $e_i$, for $i=1,2$.  These two
perpendiculars intersect at the projection of the circumcenter of
$Q$.  If the faces of $Q$ are nonobtuse,
the perpendiculars pass through the segments $e_1$ and $e_2$
respectively; and the two faces $F_1$ and $F_2$ cannot both be
negatively oriented.  \qed\enddemo

\proclaim{Lemma 2.2}  Let $Q$ be a quarter with a face $F$
along the diagonal.  Let $v$ be any
vertex not on $Q$.   If the simplex $(F,v)$
has negative orientation along $F$, then it is a quarter.
\endproclaim

A similar result holds for quasi-regular tetrahedra (Part I).

\demo{Proof}  The orientation of $F$ is determined by
the sign of the function $\chi$ (see Section I.8.2).
The face $F$ is an acute triangle, so by the explicit results for
$\chi$ in I.8.2, the function $\chi$
is increasing in the lengths of $v$ to the vertices of $F$.
  We show that $\chi\ge0$ if any of these lengths is greater than
$2.51$.  We evaluate
$$\chi(2^2,2^2,2.51^2,x^2,y^2,z^2), \quad
\chi(2^2,2.51^2,2^2,x^2,y^2,z^2),  \quad
\chi(2.51^2,2^2,2^2,x^2,y^2,z^2),$$
for $(x,y,z)\in [2,2.51]^2[2.51,\sqr8]$,
and verify that this is so.  (The minimum, which must be
attained at corner of the domain, is  $0$.)\qed\enddemo

The lemma and the results of \cite{I} imply
that if $x\in Q$ lies in the interior of Voronoi cell at a
vertex $v$ other than those of $Q$, then $v$ is a vertex of a
quarter or quasi-regular tetrahedron adjacent to $Q$.
The Voronoi cells at the vertices
of simplices in the $Q$-system cover all the simplices in the $Q$-system.

What about regions outside the $Q$-system?
A simplex $S$
in the $Q$-system may have negative orientation with respect
to a face that does not bound another simplex in the $Q$-system.
In this case, the Voronoi cell at the vertex $v_0$ opposite
this face protrudes
beyond the negatively oriented face.
More precisely, we define the tip protruding from a simplex
$S$ associated with a vertex $v_0$ of $S$ to be
the set of points that are closer to $v_0$ than to any other vertex of $S$
and that are
separated from $v_0$ by the plane through the face of $S$ opposite
$v_0$.
Each point $x$
outside the $Q$-system belongs
to finitely  many protruding tips from simplices in the $Q$-system,
say those
associated with the vertices $w_1(x),\ldots, w_k(x)$.
(Typically, this collection
of vertices is empty.)  Deleting the
vertices $w_i(x)$ from the packing, we take the Voronoi decomposition of the
remaining collection of vertices.
The point $x$ lies in the (modified) Voronoi cell at some vertex
$w(x)\ne w_i(x)$.
The set of points $x$ outside the $Q$-system such that
$v=w(x)$ will
be called the $V$-cell at $v$.  By construction, points in the
$Q$-system do not lie in any $V$-cell. Outside the $Q$-system,
$V$-cells agree with Voronoi cells except in the treatment of
protruding tips.  Occasionally, we will refer to the faces of $V$-cells
as $V$-faces to distinguish them from other types of faces,
such as those of quarters.

This is our decomposition of space: all the simplices in the $Q$-system and
all the $V$-cells.

\head{3. Scoring}\endhead

To each vertex $v$,
we attach a {\it decomposition star\/},
which is defined as the
union of the $V$-cell at $v$ with all the quasi-regular tetrahedra and
 quarters in the $Q$-system with a vertex at $v$.
Decomposition stars replace the Delaunay
stars found in earlier papers.
By construction, $V$-cells, the $Q$-system, and decomposition stars
are compatible with
standard regions.  By this, we mean in particular that the intersection
of a $V$-cell with the cone over a standard region is entirely determined
by the vertices in the cone. (See II.2.2.)
 Also, each simplex in the $Q$-system
lies over a single standard region.

A {\it standard cluster\/} attached to a standard region $P$ is the
union of the simplices in the $Q$-system over $P$ with the part of the
$V$-cell that lies over $P$.  A {\it quad cluster\/}
is the standard cluster
obtained when the standard region is a quadrilateral.

Recall that the Voronoi function $\vor(S)$ is an analytic
continuation, defined initially on simplices $S$ whose faces have
positive orientation.  Let $\sol(S)$ be the solid angle of $S$
at its distinguished vertex. Set $\doct=(\pi-4\arctan(\sqr2/5))/\sqr8$.
Set $$\vor(S) = 4(-\doct \Vol(\hat S_0) + \sol(S)/3),\tag3.1$$
where $\hat S_0\subset S$ is the Voronoi region defined in \cite{I.2}.
  An explicit formula
for $\vor(S)$ is found in \cite{I.8.6.3}.  This formula may be analytically
continued to simplices $S$ with negatively oriented faces, and $\vor(S)$
is defined in general by this analytic continuation.
Let $S_1,\ldots, S_4$ be equal to $S$ as unlabeled simplices, but with
different distinguished vertices.  Set $4\Gamma(S) = \sum_{i=1}^4\vor(S_i)$.
$\Gamma$ is called the {\it compression\/} of $S$.  The definition
here is equivalent to the one in \cite{I}.

  We define truncated versions $\vor(S,t)$ of
the Voronoi function, depending on a truncation parameter $t\le\sqr8$,
and a simplex
$S=S(y_1,\ldots,y_6)$. Set $h_i=y_i/2$, $d_i=\dih_i(S)$, the dihedral
angle along edge $i=1,2,3$.  Let $C(h,t)$ denote the compact cone
of height $h$ and circular base of area $\pi(t^2-h^2)$.
Set
$$\phi(h,t)=2(2-\doct h t (h+t))/3.$$
Then
$$2\pi(1-h/t)\phi(h,t)=(-\doct\Vol(C(h,t))+\sol(C(h,t))/3),\tag3.2$$
represents
the score of $C(h,t)$.
The solid angle of $C(h,t)$ is $2\pi(1-h/t)$, so $\phi(h,t)$
is the score per unit area.  Also, $\phi(t,t)$ is the score per unit
area of a ball of radius $t$.  That is, $\phi(t,t) = 4(-\doct\Vol/\sol + 1/3)$.

If $R=R(a,b,c)$ is a Rogers simplex (defined in I.8.6), we set
$$
\align
6\quo(R) &= (a+2c)  %
(c-a)^2\arctan(e)
        +a(b^2-a^2)e\\&-4c^3\arctan(e(b-a)/(b+c)),
\tag3.3
\endalign
$$
where $e\ge0$ is given by $e^2(b^2-a^2)=(c^2-b^2)$.
The function $\quo(R)$ (the
{\it quoin\/}
of $R$) is the volume of a wedge-like region situated above the Rogers
simplex $R$.  It is defined as
the region bounded by the four planes through the
faces of $R$ and a sphere of radius $c$ at the origin.
(See Diagram 3.4.)

\smallskip
\gram|1.7|3.4|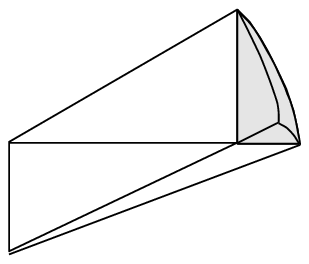|
\smallskip

We set $$
    \align
    \vor(S,t) &=
    \sol(S)\phi(t,t)
    +\sum\Sb i=1\\h_i\le t\endSb^3 d_i (1-h_i/t) (\phi(h_i,t)-
    \phi(t,t)) \\
    &-\sum_{(i,j,k)\in S_3}
    4\doct
    \quo(R(h_i,\eta(y_i,y_j,y_{k+3}),t)).
    \tag3.5
    \endalign
    $$
In the definition, we adopt the convention that $\quo(R)=0$, if
$R=R(a,b,c)$ does not exist (that is, if the condition $0<a<b<c$ is violated).
In the second sum, $S_3$ is the set of permutations on three letters.
This formula has a simple geometric interpretation
when the circumradius of $S$ is
greater than $t$ and the circumradius of each face is less than $t$.
It represents the score of the part of the Voronoi cell at the origin
that lies inside $S$ and inside a ball of radius $t$.  This can be
seen geometrically from Diagram 3.6, which depicts the intersection
of $S$ with the Voronoi cell
as three quadrilaterals forming a triangle.
The truncation in the second frame is shown as a shaded region.
The truncated volume can be decomposed
into a solid angle term, three conic terms, and six quoins (with
appropriate sign conventions).  Hence the formula for $\vor(S,t)$.

\smallskip
\gram|1.0|3.6|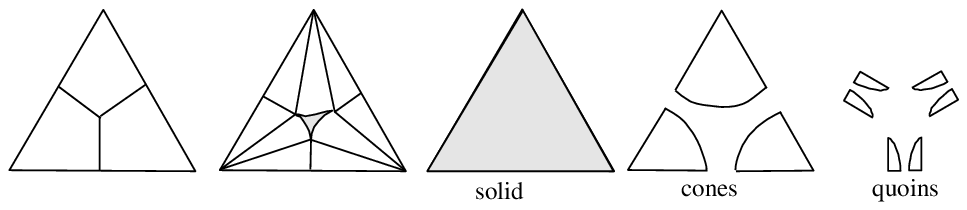|
\smallskip

Similarly, we define $\vor(P,t)$ for arbitrary standard clusters $P$.
Extending the notation in an obvious way, we have
$$
    \align
    \vor(P,t) &=
    \sol(P)\phi(t,t)
    +\sum_{|v_i|\le 2t} d_i (1-|v_i|/(2t)) (\phi(|v_i|/2,t)-
    \phi(t,t)) \\
    &-\sum_{R} 4\doct \quo(R).
    \tag3.7
    \endalign
    $$
The first sum runs over vertices in $P$ of height at most
$2t$.  The second sum runs over Rogers simplices
$R(|v_i|/2,\eta(F),t)$ in $P$, where $F=(0,v_1,v_2)$ is a
face of circumradius $\eta(F)$ at most $t$, formed by vertices
in $P$.  The constant $d_i$ is the total dihedral angle along
$(0,v_i)$ of the standard cluster.
The truncations $t=t_0=1.255=2.51/2$ and $t=\sqr2$ will
be of particular importance.
Set $A(h) = (1-h/t_0) (\phi(h,t_0)-\phi(t_0,t_0))$.

We are ready to define the
scoring of quarters and quasi-regular tetrahedra in the $Q$-system.
$\sigma(Q)$ will denote the score of a quarter or a quasi-regular
tetrahedron.  Let $S$ be a quasi-regular tetrahedron.  We set
$\sigma(S)=\Gamma(S)$ if the circumradius of $S$ is less than $1.41$,
and $\sigma(S)=\vor(S)$ otherwise.
This definition agrees with \cite{I}.

Fix a quarter $Q$.  Let $\eta^+(Q)$ be the maximum of the
circumradii of the two faces of $Q$ along the diagonal of $Q$.
Set $t_0=1.255$ and $\vor_0(Q)=\vor(Q,t_0)$.
Set
$$\mu(Q)=\cases
    \Gamma(Q),&  \text{ if }\eta^+(Q)\le\sqr2,\\
    \vor(Q),& \hbox{otherwise.}\endcases
\tag3.8
$$
If $Q$ is a flat quarter, we simply set $\sigma(Q)=\mu(Q)$.

Suppose $Q$ is upright.  Let $\hat Q$ be the upright quarter,
which is the same as $Q$ considered as an unlabeled simplex but whose
distinguished vertex lies at the opposite endpoint of the diagonal.
We say that the {\it context\/} of $Q$ is $\x(p,q)$ if there are $p-q$ quarters
along the diagonal of $Q$, and if there are $p$ anchors.
$q$ is the number of ``gaps'' between anchors around the
diagonal.
For example, the context of
a quarter in an octahedron is $\x(4,0)$.  The context
of a single quarter is $\x(2,1)$.  The only possible contexts
of upright quarters in a quad cluster are $\x(4,0)$, $\x(3,1)$, and $\x(2,1)$.
Of course, $Q$ and $\hat Q$ have the same context.
The definition of $\sigma(Q)$ depends on the context of $Q$.

context $\x(2,1)$:  Set $\sigma(Q)=\mu(Q)$.

context $\x(4,0)$:  Set $2\sigma(Q)=\mu(Q)+\mu(\hat Q)$.

other contexts:
 Set $2\sigma(Q)=\mu(Q)+\mu(\hat Q)+ \vor_0(Q) - \vor_0(\hat Q)$.

This completes the definition of $\sigma(Q)$.
Only the contexts $\x(2,1)$, $\x(3,1)$, and $\x(4,0)$
arise in the third and fifth steps of the
Kepler conjecture. (See III.2.2.)
When $\eta^+\le\sqrt2$, we say that the quarter has compression
type.  Otherwise, we say it has Voronoi type.  To say that a
quarter has compression type means that $\Gamma(Q)$ is one term of
the scoring function.  It does not mean that it is the full score.

If $Q_1$, $Q_2$, $Q_3$, and $Q_4$ are the quarter $Q$ with its
distinguished vertex placed at the four vertices of $Q$, then
it follows directly from our definitions that
$$\sum_{i=1}^4 \sigma(Q_i) = \sum_{i=1}^4\mu(Q_i)
    = \sum_{i=1}^4 \Gamma(Q_i)=4 \Gamma(Q).\tag3.9$$
Thus, the new scoring is a local reapportionment of compression,
allowing us to relate the score to the densities of packings.

Everything outside of the $Q$-system
is scored by $V$-cells.  If $P$ is a standard cluster other
than a quasi-regular tetrahedron, let $V_P$
be the intersection of the $V$-cell at the origin with
the cone over $P$.  Set
$$\vor(V_P) = 4(-\doct\Vol(V_P)+\sol(V_P)/3).\tag3.10$$
This function is not the same as the analytic Voronoi function,
defined on simplices, which is denoted in the same way.
Set $$\sigma(P) = \vor(V_P)+\sum_{Q\subset P} \sigma(Q).\tag3.11$$
The sum runs over quarters of the $Q$-system contained in $P$.
If $D^*$ is a decomposition star, set $$\sigma(D^*)=\sum_{P\subset D^*}
    \sigma(P).\tag3.12$$

Recall that the constant $\pt$,
a {\it point},  is defined as the score of a regular
quasi-regular tetrahedron with edges of length 2. We have
$\pt = 4\arctan(\sqr2/5) - \pi/3$.

\proclaim{Lemma 3.13}  A quasi-regular tetrahedron scores at most $1\,\pt$.
A quad cluster scores at most $0$, and that only for a quad cluster
whose corners have height $2$, forming a square of side $2$.
Other standard clusters have strictly negative scores.
\endproclaim

\demo{Proof}  The statement about quasi-regular tetrahedra is
found in \cite{I}.
The general context of upright quarters is established by
Calculations 3.13.3 and 3.13.4.  For the remaining quarters, it
is enough to consider $\mu(Q)$.
We claim that $\Gamma(Q)\le 0$, on quarters satisfying $\eta^+(Q)\le\sqr2$.
If the circumradius of every face of the quarter is at most $\sqr2$,
 this follows from Section II.4.5.1.
Because of this, we may assume that
the circumradius of $Q$ is greater than $\sqr2$.
The inequality $\eta^+(Q)\le\sqr2$ implies
that the faces of $Q$ along the diagonal have nonnegative orientation.
The other two faces have positive orientation, by Section I.3.4.
 Decompose the simplex into Rogers simplices as in \cite{II}
(Type IV, etc.).  The inequality $\Gamma\le0$ now
 follows from II.4 if $\eta^+\le\sqrt2$.

Assume that $\eta^+\ge\sqrt2$ and $\sigma=\vor$.
The result follows from \cite{II} if
the orientations of the sides are all positive.  In fact, we may allow
the face opposite the origin to have negative
orientation.  For the remaining cases we appeal to Calculations
3.13.1
and 3.13.2, listed in the appendix.  Calculation 3.13.1
treats flat quarters, and Calculation 3.13.2 treats the
upright quarters.

For regions outside the $Q$-system we proceed as in Part II.
We show that the score of the $V$-cell under any Voronoi face
is negative.  We adapt the fan of Part II by adding a face to the
fan if it belongs to a simplex  in the $Q$-system,
or if the circumradius of the
face is at most $\sqr2$.
Lemma II.4.4 is still valid, but its proof must
be adapted.  In the notation of \cite{II},
consider the simplex formed by $F_1$ and $F_2$.
If its circumradius is at most $\sqr2$, the argument for
small simplices in Part II
applies.  Otherwise, if the point $p$ (constructed in the Lemma)
is at most $\sqr2$ from
the vertices of $F_1$,
the face $F_2$ of the simplex has negative orientation,
giving it a circumradius greater than $\sqr2$. By
Lemma 2.2, this means that the simplex is a quarter.
If it is a quarter, since $F_2$ was included in the fan, there
is a quarter in the $Q$-system along the diagonal.  So
every quarter along the diagonal lies in the $Q$-system.  But we have
assumed that we are outside the $Q$-system.
The proof is complete.\qed\enddemo

Thus, we recover the main results of \cite{I} and \cite{II}
 under this new scoring scheme.
Set $\deltaeff(s) = 16\pi\doct/(16\pi-3s)$.  The following proposition
is a minor adaptation of Lemma I.2.1.

\proclaim{Proposition 3.14}  If every decomposition star
in a saturated packing scores at most $s<16\pi/3$, then the density
of the packing is at most $\deltaeff(s)$.  If the score of every
decomposition star is at most $8\,\pt$, then the density of
the packing is at most $\pi/\sqrt{18}$.\endproclaim

\demo{Proof}  Let $D^*(v)$ be the decomposition star around a
vertex $v$.  Let $\Lambda_N$ be the set of sphere centers
inside a large ball $B_N$ of radius $N$.
Set
$$\vor(D^*(v))=4(-\doct\Vol(V(v))+4\pi/3),$$
where $V(v)$ is the Voronoi cell around $v$.
We have
$$\sum_{\Lambda_N}\sigma(D^*(v)) = \sum_{\Lambda_N}
    \vor(D^*(v)) +O(N^2) = 4(-\doct\Vol(B_N)+
    |\Lambda_N|{4\pi\over 3}) + O(N^2).$$
This identity holds because the score of a decomposition star
$\sigma(D^*(v))$ is a local reapportionment of $\vor$.
In fact, $\Gamma$ is obtained by averaging the Voronoi volumes,
and $V$-cells are obtained from Voronoi cells
by reapportioning protruding tips
among neighboring Voronoi cells.  These modifications of the
Voronoi cells make no difference except at the boundary of $B_N$,
when we sum over $\Lambda_N$.  The term $O(N^2)$ accounts for
the boundary effects from decomposition stars that lie partially
outside $B_N$.

The inequality $\sigma(D^*(v))\le s$ gives
$$4(-\doct\Vol(B_N) + |\Lambda_N|{4\pi\over3})\le s|\Lambda_N|+O(N^2).$$
Rearranging this inequality as in the proof of Lemma I.2.1,
and taking the limit as $N$ tends to infinity,
we obtain the result.  The second statement of the
Lemma is the special case $s = 8\,\pt$.\qed\enddemo

Proposition 3.14 suggests the following conjecture.

\smallskip
\proclaim
{Conjecture 3.15}
The score of a decomposition star is at most $8\,\pt$.
\endproclaim

\proclaim{Theorem 3.16}  If a decomposition star is made
entirely of quasi-regular tetrahedra, its score is less
than $8\,\pt$.\endproclaim

\demo{Proof} Nothing has changed for quasi-regular
tetrahedra.  See \cite{I} for a proof.\qed\enddemo

\bigskip
  If the quad cluster has a diagonal of length
at most $\sqr8$ between two corners,
there are three possible decompositions.
(1) The two quarters formed by the
diagonal lie in the $Q$-system so that compression or the
Voronoi function
is used on each.  (2)  There is a second diagonal of length at most $\sqr8$,
and we use the two quarters from the second diagonal for the scoring.
(3)  There is an enclosed vertex that makes the quad cluster into
an octahedron and the four upright quarters are in the
$Q$-system.

Now suppose that neither diagonal is less than $\sqr8$ and the
quad cluster is not an octahedron.
If there is no enclosed
vertex of length at most $\sqr8$,
the quad cluster contains no quarters. An upper bound
on the score of the quad
cluster $P$ is $\vor(P,\sqr2)$.
The remaining cases are called {\it mixed\/} quad clusters.
Mixed quad clusters enclose a vertex of height at most $\sqr8$ and
do not contain flat quarters.

\bigskip
\head{4. Bounds on the Score}\endhead
\bigskip

\proclaim{Proposition 4.1}  The score of a mixed quad cluster is less
than $-1.04\,\pt$.\endproclaim

\demo{Proof}
Any enclosed vertex in a quad cluster has length
at least $2.51$ by Section III.2.2.  In particular, the anchors of an
enclosed vertex are corners of the the quad cluster.
There are no flat quarters.

We generally truncate the $V$-cell
at $\sqr2$.  This increases the score, and yet
by \cite{II} and Lemma 3.13, it breaks into pieces
whose score is nonpositive.  Thus, if we identify certain pieces
that score less than $-1.04\,\pt$, the result follows.  Nevertheless,
a few simplices will be left untruncated in the following argument.
We will leave a simplex untruncated only if we are certain that
each of its faces has positive orientation and that the simplices
sharing a face $F$ with $S$ either lie in the $Q$-system or have
positive orientation along $F$.  If these conditions hold, we may
use the Voronoi function on $S$ rather than truncation. (See
Calculations 4.1.1 and 4.1.3.)

By enclosed vertex, we now mean one of height at most $\sqrt2$.
Let $v$ be an enclosed vertex with the fewest anchors.  Consider
the part of the $V$-cell under the $V$-face determined by $v$.
If there are no anchors, under this face lies the right-circular
cone $C(h,\eta_0(h))$, where $\eta_0(h):=\eta(2h,2,2.51)$
and $|v|=2h$.  In fact, any neighboring
face corresponds to a corner of the quad cluster or to an
enclosed vertex of height at least 2.51.  In either case, the
set of points in the face's plane, at distance at most $\eta_0(h)$
from the origin, belongs to the face.
By Formula 3.2,
the score of this cone is $2\pi(1-h/\eta_0(h))\phi(h,\eta_0(h))$.
An optimization in one variable gives an upper bound of $-4.52\,\pt$,
for $1.255\le h\le \sqr2$.   This gives the bound of $-1.04\,\pt$
in this case.

If there is one anchor,  we cut the cone in half along the
plane through $(0,v)$ perpendicular to the plane containing the
anchor and $(0,v)$.  The half of the cone on the far side of the
anchor lies under the face at $v$ of the $V$-cell.  We get a bound of
$-4.52\,\pt/2 < -1.04\,\pt$.

To treat the remaining cases, we define a function $K(S)$ on certain
simplices $S$ with circumradius at least $\sqr2$.
Let $S=S(y_1,y_2,\ldots,y_6)$.  Let $R(a,b,c)$ denote a Rogers
simplex. Set
$$K(S) = K_0(y_1,y_2,y_6)+K_0(y_1,y_3,y_5)
    + \dih(S)(1-y_1/\sqr8) \phi(y_1/2,\sqr2),\tag4.2$$
where
$$\align
K_0(y_1,y_2,y_6) &= \vor(R(y_1/2,\eta(y_1,y_2,y_6),\sqr2))
    +\vor(R(y_2/2,\eta(y_1,y_2,y_6),\sqr2))\\
    &- \dih(R(y_1/2,\eta(y_1,y_2,y_6),\sqr2))
    (1-y_1/\sqr8)\phi(y_1/2,\sqr2).
\endalign$$
(If the given Rogers simplices do not exist because the
condition $0<a<b<c$ is violated,
we set the corresponding terms in these expressions to 0.)
The function $K(S)$ represents the part of the score
coming from
 the four Rogers simplices along two
of the faces of $S$, and the conic region extending out to $\sqr2$
between the two Rogers simplices along the edge $y_1$ (Diagram 4.3).

\smallskip
\gram|1.25|4.3|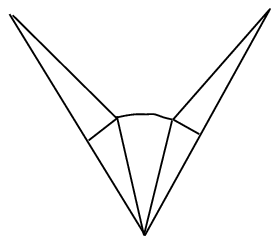|
\smallskip

Fix an enclosed vertex $v$ and draw its anchors.  Suppose that $v_1$,
a corner, is an anchor of $v$.  Assume that
the face $(0,v,v_1)$ bounds at most one upright quarter.  We
sweep around the edge $(0,v_1)$, away from the upright quarter if
there is one,  until we come to another enclosed
vertex $v'$ such that $(0,v_1,v')$ has circumradius less than $\sqr2$
or such that $v_1$ is an anchor of $(0,v')$.  If such a vertex $v'$
does not exist, we sweep all the way to $v_2$ a corner of the
quad cluster adjacent to $v_1$.

If $v'$ exists, then 4.1.1 or 4.1.2 gives the bound $-1.04\,\pt$, depending
on the size of the circumradius of $(0,v,v')$.
This allows us to assume that we do not encounter such an
enclosed vertex $v'$ whenever
we sweep away, as above, from the face formed by an anchor.

Now consider the simplex $S=(0,v_1,v_2,v)$, where $v_1$ is an
anchor of $(0,v)$.  We assume that it is not an upright quarter.
There are three alternatives.
The first is that $S$ decreases the score of the quarter
by at least $0.52\,\pt$.
This occurs if the circumradius of the face $(0,v,v_2)$ is less than
$\sqr2$ by Calculation 4.1.3, or if the circumradius of the
face is greater than $\sqr2$ by Calculation 4.1.4, provided
that the length of $(v,v_1)$ is at most $2.2$.
The second alternative is that the face $(0,v,v_1)$ of $S$ is
shared with a quarter $Q$ and that $S$ and $Q$ taken together bring
the score down by $0.52\,\pt$  (see Calculations 4.1.5 and 4.1.6).
In fact, if there are two such simplices $S$ and $S'$ along $Q$,
then the three simplices
$Q$, $S$, and $S'$ pull the score below $-1.04\,\pt$
(see Calculation 4.1.7).  The third alternative is that there is
a simplex $S'=(0,v,v,v_3)$ sharing the face $(0,v,v_1)$, which, like
$S$, scores less than $-0.31\,\pt$.  In each case, $S$
and the adjacent simplex through $(0,v,v_1)$
score less than $-0.52\,\pt$.
Since $v$ has at least two anchors, the quad cluster scores less than
$2(-0.52)\,\pt =-1.04\,\pt$.\qed\enddemo

Set $\phi_0=\phi(t_0,t_0)\approx -0.5666$.
We define $$\cro(h) = 2\pi(1-h/\eta_0(h))(\phi(h,\eta_0(h))-\phi_0).$$
It is equal to $-4\doct$ times the volume
of the region outside the sphere of radius $t_0$ and inside the
finite cone $C(h,\eta_0(h))$.  If $v$ is an enclosed vertex of height
$2h\in[2.51,\sqr8$],
 such that every other vertex $v'$ of the
standard cluster satisfies
$$\eta(|v|,|v'|,|v-v'|)\ge \eta_0(h),$$ then the
volume represented by $\cro(|v|/2)$ lies outside the truncated
$V$-cell, but inside the $V$-cell, so that if $P$ is a quad cluster,
 $$\vor(V_P) < \vor_0(V_P) + \cro(|v|/2).$$
If a vertex $v'$ satisfies $\eta(|v|,|v'|,|v-v'|)\le\eta_0(h)$, then
by the monotonicity of the circumradius of acute triangles,
$v'$ is an anchor of $v$.  This anchor clips the crown just
defined, and we add a correction term $\anc(|v'|,|v|,|v-v'|)$
to account for this.  Diagram 4.4 illustrates the terms in the
definition of $\anc()$.

\smallskip
\gram|1.7|4.4|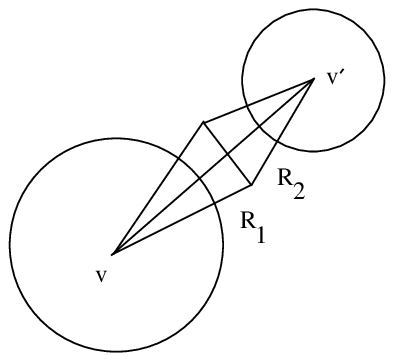|
\smallskip

Set $$\align
    \anc(y_1,y_2,y_6) &= -\dih(R_1)\cro(y_1/2)/(2\pi)
        -\sol(R_1)\phi_0+\vor(R_1)\\
    &-\dih(R_2)(1-y_2/2.51)(\phi(y_2/2,t_0)-\phi_0)
        -\sol(R_2)\phi_0 + \vor(R_2),
    \tag4.5
    \endalign
$$
where $R_i=R(y_i/2,\eta(y_1,y_2,y_6),\eta_0(y_1/2))$, for $i=1,2$.
In general, there are Rogers simplices on both sides of the
face $(0,v,v')$, and this gives a factor of 2.
For example, if $v$ has a single anchor $v'$, then
$$\vor(V_P) < \vor_0(V_P) + \cro(|v|/2) + 2\anc(|v|,|v'|,|v-v'|).$$
However,
if the anchor gives a face of an upright quarter,
only one side of the face lies in the $V$-cell,
 so that the factor of 2 is not required.
For example,
$v'$ has context $\x(2,1)$ with upright quarter $Q$, and if there
are no other enclosed vertices, and if
$v',v''$ are the anchors along the faces of the
quarter,  then
$$\align\vor(V_P)&< \vor_0(V_P) +(1-\dih(Q)/(2\pi))\cro(|v|/2)\\
    &+\anc(|v|,|v'|,|v-v'|)+\anc(|v|,|v''|,|v-v''|).
    \endalign$$
In general, when there are multiple anchors around the same enclosed vertex
$v$, we add a term $(2-k)\anc$ for each anchor, where
$k\in\{0,1,2\}$ is the number of
quarters bounded by the face formed by the anchor.
We must be cautious in the use of this formula.
If the circumradius
of $(0,v,v',v'')$ is less than $\eta_0(|v|/2)$, the
Rogers simplices used to define the terms $\anc()$ at $v'$ and $v''$
overlap.  When this occurs, the geometric decomposition on which
the correction terms $\anc()$ are based is no longer valid.  In this
case, other methods must be used.

\smallskip
\gram|2|4.6|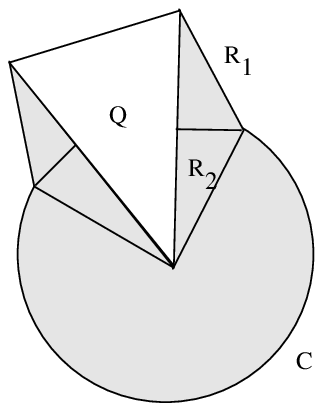|
\smallskip

\smallskip
If $P$ is a mixed quad cluster, let $P_0$ be the new quad cluster obtained
by removing all the enclosed vertices.  We define a $V$-cell $V_{P_0}$
of $P_0$ and the truncation of $V_{P_0}$ at $t_0$. We take its score
$\vor_0(P_0)$  as we do for standard clusters.  $P_0$ does
not contain any quarters.

\proclaim{Proposition 4.7}  If $P$ is a mixed quad cluster,
$\sigma(P) < \vor_0(P_0)$.
\endproclaim

\bigskip
\demo{Proof}  Suppose there exists an enclosed vertex that
has context $\x(2,1)$; that is, there is a single upright quarter
$Q=S(y_1,y_2,\ldots,y_6)$
and no additional anchors.  In this context $\sigma(Q)=\mu(Q)$.
Let $v$ be the enclosed vertex.  To compare $\sigma(P)$ and $\vor_0(P_0)$,
consider the $V$-cell near $Q$.
The quarter $Q$ cuts a wedge of angle $\dih(Q)$ from the crown at $v$.
There is an anchor term for the two anchors of $v$ along the faces
of $Q$.  Let $V_P^v$ be the truncation at height $t_0$ of $V_P$ under the
$V$-face determined by $v$ and under the four Rogers simplices stemming
from the two anchors.  (Diagram 4.6 shades the truncated parts of the
quad cluster.) As a consequence
\smallskip
$$\vor(V_P) <(1-\dih(Q)/(2\pi))\cro(y_1/2)+\anc(y_1,y_2,y_6)
+\anc(y_1,y_3,y_5) +\vor(V_P^v).\tag4.8$$
Combining this inequality with Calculations 4.7.2, 4.7.3, and 4.7.4,
    we find
$$\vor(V_P) +\mu(Q) < \vor(V_P^v) +\vor_0(Q).\tag4.9$$

Now suppose there is an enclosed vertex $v$ with context $\x(3,1)$.
Let the quad cluster have corners $v_1$, $v_2$, $v_3$, $v_4$,
ordered consecutively.  Suppose the two quarters along $v$ are
$Q_1=(0,v,v_1,v_2)$ and $Q_2=(0,v,v_2,v_3)$.  We consider two
cases.

\noindent
Case 1:  $\dih(Q_1)+\dih(Q_2)<\pi$ or $\rad(0,v,v_1,v_3)>\eta(|v|,2,2.51)$.
In this case, the use of correction terms to the crown are
legitimate (in relation to the note of caution about
the possible overlap of Rogers simplices).  Proceeding as
in context $\x(2,1)$, we find that
\smallskip
$$\vor(V_P) < (1-(\dih(Q_1)+\dih(Q_2))/(2\pi))\cro(|v|/2)
    +\anc(F_1) +\anc(F_2) +\vor(V_P^v).\tag4.10$$
Here $V_P^v$ is defined by the truncation at height $t_0$
 under the $V$-face determined by $v$ and
under the Rogers simplices stemming from the side of $F_i$ that
occur in the definition of $\anc$.
Also, $\anc(F_i)=\anc(y_i,y_j,y_k)$ for a face $F_i$ with edges $y_i$
along an upright quarter.
   By Calculation 4.7.1 applied to both $Q_1$ and $Q_2$, we have
$$\vor(V_P) +\sum_{i=1}^2\sigma(Q_i)
    < \vor(V_P^v) + \sum_{i=1}^2 \vor_0(Q_i).\tag4.11$$
That is, by truncating near $v$, and changing the scoring
of the quarters to $\vor_0$, we obtain an upper bound on the score.

\noindent
Case 2:  $\dih(Q_1)+\dih(Q_2)\ge\pi$ and
    $\rad(0,v,v_1,v_3)\le \eta_0(|v|/2)$.
The anchor terms cannot be used here.
In the mixed case, $\sqr8<|v_1-v_3|$, so
$$\sqr2<{1\over2}|v_1-v_3|\le\rad \le \eta_0(|v|/2),$$
and this implies $|v|\ge 2.696$.
We have $$\sum_{i=1}^2 \sigma(Q_i) <
\sum_{i=1}^2 \vor_0(Q_i) + \sum_{i=1}^2 0.01(\pi/2-\dih(Q_i))< \sum_{i=1}^2
\vor_0(Q_i)$$
by Calculation 4.7.5.
Inequality 4.11 holds, for $V_P^v=V_P$.

In the general case, we run over all enclosed vertices $v$
and truncate around each vertex.  For each vertex we obtain
4.9 or 4.11.   These inequalities can be coherently combined
over multiple enclosed vertices because the $V$-faces were associated
with different vertices $v$
and none of the Rogers simplices used in the terms $\anc()$ overlap.
More precisely, if $Z$ is a set of enclosed vertices, set
$V_P^Z = \cap_{v\in Z} V_P^v$, and $V_P^{v,Z} = V_P^Z\cap V_P^v$.
Coherence means that we obtain valid inequalities by adding the superscript
$Z$ to $V_P$ and $V_P^v$ in Inequalities 4.9 and 4.11, if $v\not\in Z$.
In sum,  $\sigma(P) <\vor_0(P_0)$.\qed\enddemo

\newpage

\bigskip
\head{Appendix 1.}\endhead

The following inequalities have been proved by computer using
interval methods.  The standard methods described in \cite{I} have
been used, together with various improvements in method that will
be described elsewhere.
  Let $S=S(y)=S(y_1,\ldots,y_6)$.
Set $\eta_{234}=\eta(y_2,y_3,y_4)$,
and $\eta_{126} = \eta(y_1,y_2,y_6)$.
The function $K(S)$ is introduced in
Section 4.1.
\parindent=0pt\hbox{}

\smallskip

{\bf Calculation 3.13.1.}
$\vor(S)\le 0$, for $y\in[2,2.51]^3[2.51,\sqr8][2,2.51]^2$
if the orientation is negative for the
face containing the origin and the long edge.

{\bf Calculation 3.13.2.}
$\vor(S)\le 0$, for $y\in[2.51,\sqr8][2,2.51]^5$.

{\bf Calculation 3.13.3.}
$2\Gamma(S) + \vor_0(S)-\vor_0(\hat S) \le0$,
    for all upright quarters $S$.

{\bf Calculation 3.13.4.}
$\vor(S) + \vor(\hat S)+ \vor_0(S)-\vor_0(\hat S) \le0$,
    for all upright quarters $S$.

{\bf Calculation 4.1.1.}
$\vor(S) < -1.04\,pt$,
    provided $\eta_{234}\le\sqr2$
    and $$y\in[2,2.51][2.51,2.7]^2[2,2.32][2,2.51]^2,$$
    or provided $\eta_{234},\eta_{126}\le\sqr2$
    and $$y\in[2,2.51][2.51,2.7]^2[2,2.32][2,2.51][2.51,2.7].$$
(We have $y_4\le 2.32$ because otherwise
$\eta_{234}>\eta(2.51,2.51,2.32)>\sqr2$.  Similarly, $y_2,y_3\le 2.7$;
otherwise $\eta_{234}>\eta(2.51,2.7,2)>\sqr2$. Similarly, $y_6\le 2.7$
if $\eta_{126}\le\sqr2$.)

{\bf Calculation 4.1.2.}
$K(S) < -1.04\,\pt$,
    provided $\eta_{234}\ge\sqr2$
    and $$y\in[2,2.51][2.51,\sqr8]^2[2,2.51]^3,$$
    or provided $\eta_{234}\ge\sqr2\ge\eta_{126}$
    and $$y\in[2,2.51][2.51,\sqr8]^2[2,2.51]^2[2.51,2.7].$$

{\bf Calculation 4.1.3.}  $\vor(S) < -0.52\,\pt$,
    provided $\eta_{234}\le\sqr2$
    and $$y\in[2,2.51]^2[2.51,2.7]^2[2,2.51]^2.$$
(We have $y_4\ge 2.51$, because $S$ is assumed not to be a quarter.)

{\bf Calculation 4.1.4.}  $K(S) < -0.52\,\pt$, provided
    $\eta_{234}\ge\sqr2$
    and $$y\in[2,2.51]^2[2.51,\sqr8]^2[2,2.2][2,2.51].$$

{\bf Calculation 4.1.5.}  $K(S) < -0.31\,\pt$, provided
    $\eta_{234}\ge\sqr2$
    and $$y\in[2,2.51]^2[2.51,\sqr8]^2[2,2.51]^2.$$

{\bf Calculation 4.1.6.}  $\sigma(Q) < -0.21\,\pt$
    in the contexts $\x(2,1)$ and $\x(3,1)$
    for
    $$y\in[2.51,\sqr8][2,2.51]^4[2.2,2.51].$$

{\bf Calculation 4.1.7.}  $\sigma(Q) < -0.42\,\pt$
    in the context $\x(2,1)$
    for $$y\in[2.51,\sqr8][2,2.51]^3[2.2,2.51]^2.$$

{\bf Calculation 4.7.1.}
    For all upright quarters $Q$,
    $$\mu(Q)+\mu(\hat Q)
    +(1-\dih(Q)/\pi)\cro(y_1/2)+2\anc(y_1,y_2,y_6)
        < \vor_0(Q)+\vor_0(\hat Q).$$

{\bf Calculation 4.7.2.}
    $\cro(h) < -0.1378$, for $h\in [1.255,\sqr2]$.

{\bf Calculation 4.7.3.}
    $\anc(y_1,y_2,y_6) < 0.0263$, for $$(y_1,y_2,y_6)
        \in [2.51,\sqr8][2,2.51]^2.$$

{\bf Calculation 4.7.4.}
    $\mu(Q) + (1-\dih(Q)/2\pi)(-0.1378) + 2(0.0263) < \vor_0(Q)$,
        for all upright quarters $Q$.

{\bf Calculation 4.7.5.}
    $\mu(Q)+\mu(\hat Q) <
    \vor_0(Q)+\vor_0(\hat Q) +
    0.02 (\pi/2-\dih(Q))$, for $y\in [2.69,\sqr8][2,2.51]^5$.

\bigskip

\newpage

\bigskip
\head{Appendix 2. Compatibility Notes}\endhead

\parskip=\baselineskip
It has been useful to make various changes in the program that was
published in {\it Sphere Packings I}. This appendix makes a few
comments about the global compatibility of the results and
terminology from various papers.

The definition of {\it quasi-regular octahedron} in Sphere
Packings I is obsolete.  The definition that is used appears in
Section I of this paper.  Also, there is an old definition of {\it
standard cluster\/} for Delaunay stars that should be replaced
with the standard cluster in a decomposition star in \cite{F}.

The Sections I.8.6.4, I.8.6.5, I.8.6.6, I.8.6.7 are no longer
needed because of improvements in the numerical methods used to
calculate the Voronoi function.  Also, Lemma I.9.1.1 can now be
verified quickly by computer, so the technical proof that is given
is no longer needed.  Lemmas I.9.17 and I.9.18 are proved by a
long argument that is no longer necessary because of improvements
in numerical methods.

Many of the papers rely on the $\arctan$ formula for the dihedral
angle, rather than the $\arccos$ formula that appears in I.8.
$$\dih(S) = \pi/2 +\arctan(-\Delta_4/\sqrt{4x_1\Delta}).$$
This leads to simple formulas for the derivatives of the dihedral
angles that have been used extensively throughout the collection
without explicit mention
$$\partial_2\dih = -y_1\Delta_3/(u_{126}\sqrt{\Delta}),$$
etc.

In {\it Sphere Packings II\/} the notion of a {\it small} simplex
is made obsolete by the constructions of \cite{F}.  Delaunay stars
are replaced by decomposition stars.  {\it Restricted cells\/} are
replaced with $V$-cells in \cite{F}.  Simplices of {\it
compression type\/} undergo a small change in meaning when the
scoring functions are adjusted in \cite{F}. In \cite{II}, in
constructing the standard regions, we remove all arcs that do not
bound a region, but in the classification of standard regions in a
later paper these arcs will not be removed.

Lemma II.2.2 can be proved by simpler means.  After the first
paragraph of the proof, we observe that $S=(v_0,v_1,v_2,w)$ has
negative orientation along $F=(v_0,v_1,v_2)$.  Hence $S$ is a
quasi-regular tetrahedron by I.3.4.  Various lemmas are revised in
\cite{F} to account for the change in decomposition. (Lemma
II.2.4, Section II.3.1, Lemma II.3.2, Theorem II.4.1.b). Several
of the cases in II.4.5.2 are unnecessary in light of the revisions
in \cite{F}.  The technical results in the appendix can now be
obtained quickly by computer.

When we say that a simplex has compression type, it means the the
scoring rule for $\eta^+(Q)\le\sqrt2$ is used.  Here $\eta^+$ is
the function of Section 3.  To say that a simplex has compression
type implies that the compression function is one term of the
scoring function.  But there will often be various correction
terms, so that the scoring function need not be identical with the
compression function.  Similar comments apply to simplices of
Voronoi type, which means precisely that $\eta^+(Q)>2\sqrt2$. It
means loosely that the Voronoi function appears as part of the
scoring function. In general, the terms Voronoi, Voronoi scoring,
Voronoi function, and so forth are used loosely for objects
related to the $V$-cells in the decomposition star.

\newpage
\head{References.}\endhead

[FT] L. Fejes T\'oth, Lagerungen in der Ebene auf der Kugel und im Raum,
    Berlin, Springer-Verlag, 1953.

[V] S.P. Ferguson, Sphere Packings, V, thesis, University of Michigan, 1997.

[H1] T.C. Hales, Remarks on the density of sphere packings in three
    dimensions, Combinatorica, 13(2), 1993, 181-197.

[H2] T.C. Hales, the Sphere Packing Problem, J. Comp. App. Math. 44, 1992,
    41-76.

[I]  T.C. Hales, Sphere Packings I, Disc.
        Comp. Geom. 1997, 17:1-51.  [II] Sphere Packings II,
    DCG. 1997, 18:135-149. [III] Sphere Packings III, preprint.

\enddocument
\bye